\documentclass[12pt]{article}%
\usepackage{amsmath}
\usepackage{amsfonts}
\usepackage{amssymb}
\usepackage{graphicx}%
\newtheorem{theorem}{Theorem}

\newtheorem{conclusion}[theorem]{Conclusion}

\newtheorem{corollary}[theorem]{Corollary}

\newtheorem{definition}[theorem]{Definition}
\newtheorem{example}[theorem]{Example}

\newtheorem{remark}[theorem]{Remark}

\newenvironment{proof}[1][Proof]{\noindent\textbf{#1.} }{\ \rule{0.5em}{0.5em}}
\begin{document}

\title{An extension of the classical derivative}
\author{Diego Dominici \thanks{e-mail: dominicd@newpaltz.edu}\\Department of Mathematics\\State University of New York at New Paltz\\75 S. Manheim Blvd. Suite 9\\New Paltz, NY 12561-2443\\USA\\Phone: (845) 257-2607\\Fax: (845) 257-3571 }
\date{}
\maketitle

\section{Introduction}

Undergraduate students attending Calculus I in their first semester at college
(or increasingly in high school \cite{Bressoud}) will learn how to compute the
derivative of a function of a real variable, by means of the formula%
\begin{equation}
f^{\prime}(c)=\underset{x\rightarrow c}{\lim}\frac{f(x)-f(c)}{x-c}%
.\label{classder}%
\end{equation}
Some of them may be surprised by the fact that while the derivative of
functions like $x$ or $x^{2}$ is defined at $x=0,$ a well behaved function
like $\sqrt{x}$ is not differentiable at $x=0.$

Advancing to Calculus II, they will then learn how to find a Taylor or
Maclaurin series to represent an elementary function. Once again, they will
find that everything seems to break down when they try to expand functions
such as $\sin\left(  \sqrt{x}\right)  $ at $x=0.$ Of course, a (bad!) way of
getting around this is to compute the Maclaurin series of $\sin(x)$%
\[
\sin(x)=%
%TCIMACRO{\dsum \limits_{n=0}^{\infty}}%
%BeginExpansion
{\displaystyle\sum\limits_{n=0}^{\infty}}
%EndExpansion
\frac{\left(  -1\right)  ^{n}}{\left(  2n+1\right)  !}x^{2n+1},
\]
and then replace $x$ by $\sqrt{x},$ to obtain%
\begin{equation}
\sin(\sqrt{x})=%
%TCIMACRO{\dsum \limits_{n=0}^{\infty}}%
%BeginExpansion
{\displaystyle\sum\limits_{n=0}^{\infty}}
%EndExpansion
\frac{\left(  -1\right)  ^{n}}{\left(  2n+1\right)  !}x^{n+\frac{1}{2}},\quad
x\geq0.\label{taylsin}%
\end{equation}
Inquisitive students may ask themselves, What is wrong with fractional powers?

A few of them will finally understand when they take a course in Complex
Analysis and learn about branch points, Riemann surfaces, etc. Unfortunately,
those will constitute a very small percentage of the whole population of
Calculus students.

Over the years, there have been many generalizations of the classical
derivative \cite{MR0197632}. Most of them are beyond the reach of students at
the Calculus I level. In this paper, we propose an extension of the classical
derivative that would be easy to teach to Calculus I students and which
provides a simple way of computing the expansion of functions at branch points.

In the remainder of this paper, we will always refer to functions whose
derivative defined by (\ref{classder}) exists, as differentiable functions.

\section{Definition and examples}

\begin{definition}
Let $I\subset\mathbb{R}$ be an interval, $f:I\rightarrow\mathbb{R}$ and $c\in
I.$ Let $\alpha>0$ and $\omega$ be a real number. We say that a number
$D_{\alpha,\omega}\left[  f\right]  \left(  c\right)  $ is the $\left(
\alpha,\omega\right)  $\textbf{-derivative} of $f$ at $c$ if the limit%
\begin{equation}
D_{\alpha,\omega}\left[  f\right]  \left(  c\right)  =\underset{x\rightarrow
c}{\lim}\frac{f(x)-f(c)}{\left(  x-\omega\right)  ^{\alpha}-\left(
c-\omega\right)  ^{\alpha}} \label{def1}%
\end{equation}
exists. In this case, we say that $f$ is $\left(  \alpha,\omega\right)
$\textbf{-differentiable} at $c.$

We could also define $D_{\alpha,\omega}\left[  f\right]  \left(  c\right)  $
by the equivalent limit%
\[
D_{\alpha,\omega}\left[  f\right]  \left(  c\right)  =\underset{h\rightarrow
0}{\lim}\frac{f(c+h)-f(c)}{\left(  c+h-\omega\right)  ^{\alpha}-\left(
c-\omega\right)  ^{\alpha}}.
\]

\end{definition}

We will use the notation $D_{\alpha,\omega}\left[  f\right]  $ to denote the
function whose domain is a subset of $I$ and whose value at $c$ is
$D_{\alpha,\omega}\left[  f\right]  \left(  c\right)  .$ Higher derivatives
will be denoted by $D_{\alpha,\omega}^{k}\left[  f\right]  \left(  c\right)  $
and defined by $D_{\alpha,\omega}^{0}\left[  f\right]  \left(  c\right)
=f(c)$ and%
\[
D_{\alpha,\omega}^{k+1}\left[  f\right]  \left(  c\right)  =D_{\alpha,\omega
}\left[  D_{\alpha,\omega}^{k}\left[  f\right]  \right]  \left(  c\right)
,\quad k=1,2,\ldots.
\]

\begin{remark}
It is clear from the definition (\ref{def1}) that if $f$ is differentiable at
$x=c,$ then%
\[
D_{1,\omega}\left[  f\right]  \left(  c\right)  =\frac{df}{dx}\left(
c\right)
\]
for all $\omega.$
\end{remark}

\begin{example}
\begin{enumerate}
\item The derivative of the function $f(x)=\left(  x-\omega\right)  ^{\alpha}$
at $x=\omega$ is undefined for all $0<\alpha<1.$ Could $\left(  \alpha
,\omega\right)  $-derivatives do better? Since
\[
\frac{f(x)-f(c)}{\left(  x-\omega\right)  ^{\alpha}-\left(  c-\omega\right)
^{\alpha}}=1,
\]
we conclude from (\ref{def1}) that
\begin{equation}
D_{\alpha,\omega}\left[  \left(  x-\omega\right)  ^{\alpha}\right]  \left(
c\right)  =1.\label{basic}%
\end{equation}
This shows that the function $f(x)=\left(  x-\omega\right)  ^{\alpha}$ plays
the same fundamental role as $f(x)=x$ does for the classical derivative. In
particular, we have
\[
D_{\frac{1}{2},0}\left[  \sqrt{x}\right]  =1.
\]

\item At $x=0,$ none of the derivatives of the function $f(x)=\sin\left(
\sqrt{x}\right)  $ exist. On the other hand,%
\[
D_{\frac{1}{2},0}\left[  f\right]  \left(  0\right)  =\underset{x\rightarrow
0}{\lim}\frac{\sin\left(  \sqrt{x}\right)  }{\sqrt{x}}=1.
\]
We will prove in the next section that
\begin{equation}
D_{\frac{1}{2},0}^{k}\left[  f\right]  \left(  0\right)  =\left\{
\begin{array}
[c]{c}%
0,\quad\quad\quad\quad k\text{ even}\\
\left(  -1\right)  ^{\frac{k-1}{2}},\quad k\text{ odd}%
\end{array}
\right.  ,\quad k=1,2,\ldots,\label{sin}%
\end{equation}
and therefore we can write%
\[
\sin\left(  \sqrt{x}\right)  =%
%TCIMACRO{\dsum \limits_{n=0}^{\infty}}%
%BeginExpansion
{\displaystyle\sum\limits_{n=0}^{\infty}}
%EndExpansion
D_{\frac{1}{2},0}^{n}\left[  f\right]  \left(  0\right)  \frac{x^{\frac{n}{2}%
}}{n!},\quad x\geq0
\]
rather than (\ref{taylsin}).
\end{enumerate}
\end{example}

\section{Properties}

In this section we show some of the basic properties of the $\left(
\alpha,\omega\right)  $-derivative. We start with a theorem that generalizes
the usual proof that continuity is a necessary condition for the existence of
a derivative and the proofs of the basic differential rules.

\begin{theorem}
Let $I\subset\mathbb{R}$ be an interval, $f,g,\rho:I\rightarrow\mathbb{R}$ and
$c\in I.$ Let $\mathfrak{D}\left[  f\right]  \left(  c\right)  $ be defined by%
\[
\mathfrak{D}\left[  f\right]  \left(  c\right)  =\underset{x\rightarrow
c}{\lim}\frac{f(x)-f(c)}{\rho(x)},
\]
where $\rho(x)$ is continuous at $x=c,$ and $\rho(c)=0.$ Suppose that
$\mathfrak{D}\left[  f\right]  \left(  c\right)  $ and $\mathfrak{D}\left[
g\right]  \left(  c\right)  $ exist. Then,

\begin{enumerate}
\item[(i)] $f$ is continuous at $x=c.$

\item[(ii)] $\mathfrak{D}$ is linear, i.e.,%
\[
\mathfrak{D}\left[  af+bg\right]  (c)=a\mathfrak{D}\left[  f\right]  \left(
c\right)  +b\mathfrak{D}\left[  g\right]  \left(  c\right)  .
\]

\item[(iii)] Product rule:%
\[
\mathfrak{D}\left[  f\times g\right]  (c)=\mathfrak{D}\left[  f\right]
\left(  c\right)  \times g(c)+f(c)\times\mathfrak{D}\left[  g\right]  \left(
c\right)  .
\]

\item[(iv)] Quotient rule:%
\[
\mathfrak{D}\left[  \frac{f}{g}\right]  (c)=\frac{\mathfrak{D}\left[
f\right]  \left(  c\right)  \times g(c)-f(c)\times\mathfrak{D}\left[
g\right]  \left(  c\right)  }{\left[  g(c)\right]  ^{2}},
\]
as long as $g(c)\neq0.$
\end{enumerate}
\end{theorem}

\begin{proof}

\begin{enumerate}
\item[(i)] Since
\[
f(x)=f(c)+\frac{f(x)-f(c)}{\rho(x)}\rho(x),\quad x\neq c,
\]
we have%
\[
\underset{x\rightarrow c}{\lim}f(x)=f(c)+\mathfrak{D}\left[  f\right]  \left(
c\right)  \underset{x\rightarrow c}{\lim}\rho(x)=f(c).
\]

\item[(ii)] By definition,%
\begin{align*}
\mathfrak{D}\left[  af+bg\right]  (c)  &  =\underset{x\rightarrow c}{\lim
}\frac{\left(  af+bg\right)  (x)-\left(  af+bg\right)  (c)}{\rho(x)}\\
&  =\underset{x\rightarrow c}{\lim}a\frac{f(x)-f(c)}{\rho(x)}+b\frac
{g(x)-g(c)}{\rho(x)}\\
&  =a\mathfrak{D}\left[  f\right]  \left(  c\right)  +b\mathfrak{D}\left[
g\right]  \left(  c\right)  ,
\end{align*}
by the linearity of the limit operation.

\item[(iii)] We have
\begin{align*}
\mathfrak{D}\left[  f\times g\right]  (c)  &  =\underset{x\rightarrow c}{\lim
}\frac{\left(  f\times g\right)  (x)-\left(  f\times g\right)  (c)}{\rho(x)}\\
&  =\underset{x\rightarrow c}{\lim}\frac{f(x)g(x)-f(c)g(x)}{\rho(x)}%
+\frac{f(c)g(x)-f(c)g(c)}{\rho(x)}\\
&  =\underset{x\rightarrow c}{\lim}\frac{f(x)-f(c)}{\rho(x)}g(x)+f(c)\frac
{g(x)-g(c)}{\rho(x)}.
\end{align*}
>From part (i), we know that $\underset{x\rightarrow c}{\lim}g(x)=g(c).$ Thus,
\begin{align*}
&  \underset{x\rightarrow c}{\lim}\frac{f(x)-f(c)}{\rho(x)}g(x)+f(c)\frac
{g(x)-g(c)}{\rho(x)}\\
&  =\mathfrak{D}\left[  f\right]  \left(  c\right)  \times g(c)+f(c)\times
\mathfrak{D}\left[  g\right]  \left(  c\right)  .
\end{align*}

\item[(iv)] In this case,%
\begin{align*}
\mathfrak{D}\left[  \frac{f}{g}\right]  (c)  &  =\underset{x\rightarrow
c}{\lim}\frac{\left(  \frac{f}{g}\right)  (x)-\left(  \frac{f}{g}\right)
(c)}{\rho(x)}=\underset{x\rightarrow c}{\lim}\frac{f(x)g(c)-g(x)f(c)}%
{g(c)g(x)\rho(x)}\\
&  =\underset{x\rightarrow c}{\lim}\frac{1}{g(c)g(x)}\frac
{f(x)g(c)-f(c)g(c)+f(c)g(c)-g(x)f(c)}{\rho(x)}\\
&  =\underset{x\rightarrow c}{\lim}\frac{1}{g(c)g(x)}\left[  \frac
{f(x)-f(c)}{\rho(x)}g(c)-f(c)\frac{g(x)-g(c)}{\rho(x)}\right] \\
&  =\frac{1}{\left[  g(c)\right]  ^{2}}\left\{  \mathfrak{D}\left[  f\right]
\left(  c\right)  \times g(c)-f(c)\times\mathfrak{D}\left[  g\right]  \left(
c\right)  \right\}  .
\end{align*}

\end{enumerate}
\end{proof}

\begin{remark}
Taking
\[
\rho(x)=\left(  x-\omega\right)  ^{\alpha}-\left(  c-\omega\right)  ^{\alpha}%
\]
in the previous theorem, we find that all the results proven hold for the
$\left(  \alpha,\omega\right)  $-derivative.
\end{remark}

We shall now find a relation between $\left(  \alpha,\omega\right)
$-derivatives with different values of $\alpha$ and $\omega.$ This will
provide us with a useful way of computing $\left(  \alpha,\omega\right)
$-derivatives of differentiable functions.

\begin{theorem}

\begin{description}
\item[(i)] Let $c\neq\omega,\zeta.$ Then, $f$ is $\left(  \alpha
,\omega\right)  $-differentiable at $c$ if and only $f$ is $\left(
\beta,\zeta\right)  $-differentiable at $c.$ In this case,%
\[
D_{\alpha,\omega}\left[  f\right]  \left(  c\right)  =\frac{\beta}{\alpha
}\frac{\left(  c-\zeta\right)  ^{\beta-1}}{\left(  c-\omega\right)
^{\alpha-1}}D_{\beta,\zeta}\left[  f\right]  \left(  c\right)  .
\]

\item[(ii)] If $f$ is $\left(  \beta,\omega\right)  $-differentiable at
$\omega$ then,%
\[
D_{\alpha,\omega}\left[  f\right]  \left(  \omega\right)  =\left\{
\begin{array}
[c]{c}%
0\quad\quad\quad\quad\quad\quad\text{if }\alpha<\beta\\
\text{undefined\quad if }\alpha>\beta
\end{array}
\right.  .
\]

\end{description}
\end{theorem}

\begin{proof}

\begin{description}
\item[(i)] It follows from (\ref{def1}) that%
\begin{align*}
D_{\alpha,\omega}\left[  f\right]  \left(  c\right)   &  =\underset
{x\rightarrow c}{\lim}\frac{f(x)-f(c)}{\left(  x-\omega\right)  ^{\alpha
}-\left(  c-\omega\right)  ^{\alpha}}\\
&  =\underset{x\rightarrow c}{\lim}\frac{f(x)-f(c)}{\left(  x-\zeta\right)
^{\beta}-\left(  c-\zeta\right)  ^{\beta}}\frac{\left(  x-\zeta\right)
^{\beta}-\left(  c-\zeta\right)  ^{\beta}}{\left(  x-\omega\right)  ^{\alpha
}-\left(  c-\omega\right)  ^{\alpha}}.
\end{align*}
>From L'Hopital's rule,%
\[
\underset{x\rightarrow c}{\lim}\frac{\left(  x-\zeta\right)  ^{\beta}-\left(
c-\zeta\right)  ^{\beta}}{\left(  x-\omega\right)  ^{\alpha}-\left(
c-\omega\right)  ^{\alpha}}=\underset{x\rightarrow c}{\lim}\frac{\beta\left(
x-\zeta\right)  ^{\beta-1}}{\alpha\left(  x-\omega\right)  ^{\alpha-1}},
\]
and the result follows immediately.

\item[(ii)] We have
\begin{align*}
D_{\alpha,\omega}\left[  f\right]  \left(  \omega\right)   &  =\underset
{x\rightarrow\omega}{\lim}\frac{f(x)-f(\omega)}{\left(  x-\omega\right)
^{\alpha}}=\underset{x\rightarrow\omega}{\lim}\frac{f(x)-f(\omega)}{\left(
x-\omega\right)  ^{\beta}}\left(  x-\omega\right)  ^{\beta-\alpha}\\
&  =D_{\beta,\omega}\left[  f\right]  \left(  \omega\right)  \times0,
\end{align*}
as long as $\alpha<\beta.$
\end{description}
\end{proof}

When $\beta=1,$ we obtain the following result.

\begin{corollary}
If $f$ is differentiable at $x,$ then%
\begin{equation}
D_{\alpha,\omega}\left[  f\right]  \left(  x\right)  =\frac{1}{\alpha}\left(
x-\omega\right)  ^{1-\alpha}f^{\prime}\left(  x\right)  \label{deriv}%
\end{equation}
for all $x\neq\omega.$
\end{corollary}

\begin{example}

\begin{enumerate}
\item Let $f(x)=\left(  x-a\right)  ^{p}.$ Then, we obtain from (\ref{deriv})%
\[
D_{\alpha,\omega}\left[  f\right]  \left(  x\right)  =\frac{p}{\alpha}\left(
x-\omega\right)  ^{1-\alpha}\left(  x-a\right)  ^{p-1}.
\]
If we take $p=\alpha$ and $a=\omega,$ we recover (\ref{basic}).

\item Let's apply (\ref{deriv}) to the function $f(x)=\sin\left(  \sqrt
{x}\right)  .$ We have%
\[
D_{\alpha,\omega}\left[  f\right]  \left(  x\right)  =\frac{1}{\alpha}\left(
x-\omega\right)  ^{1-\alpha}\frac{\cos\left(  \sqrt{x}\right)  }{2\sqrt{x}}.
\]
In particular,%
\[
D_{\frac{1}{2},0}\left[  f\right]  \left(  x\right)  =\cos\left(  \sqrt
{x}\right)  ,
\]
and using (\ref{deriv}) again, we get%
\[
D_{\frac{1}{2},0}^{2}\left[  f\right]  \left(  x\right)  =-\sin\left(
\sqrt{x}\right)  .
\]
In a similar fashion, we find
\[
D_{\frac{1}{2},0}^{k}\left[  f\right]  \left(  x\right)  =\left\{
\begin{array}
[c]{c}%
\left(  -1\right)  ^{\frac{k}{2}}\sin\left(  \sqrt{x}\right)  ,\quad k\text{
even}\\
\left(  -1\right)  ^{\frac{k-1}{2}}\cos\left(  \sqrt{x}\right)  ,\quad k\text{
odd}%
\end{array}
\right.  ,
\]
which proves (\ref{sin}).
\end{enumerate}
\end{example}

\begin{conclusion}
\strut We have extended the usual definition of a derivative in a way that
Calculus I students can easily comprehend and which allows calculations at
branch points. In a forthcoming paper, we will prove the $\left(
\alpha,\omega\right)  $-versions of the classical theorems (Rolle, Lagrange,
Cauchy, L'Hopital, Taylor, etc.) and consider possible extensions to complex variables.
\end{conclusion}

\end{document}